\documentclass[a4paper,11pt,sans]{article}

\usepackage[top=2cm,bottom=2cm,left=2cm,right=2cm]{geometry}
\usepackage[utf8]{inputenc}
\usepackage[english]{babel}
\usepackage{graphicx}
\usepackage{bm}
\usepackage{amsmath}
\usepackage{epigraph}
\usepackage{amsthm}
\usepackage{amssymb}
\usepackage{authblk}
\theoremstyle{definition}

\title{Theoretical analysis of a derivative free control based continuation algorithm with path following capability for autonomous systems.}
\author[1]{Etienne Gourc}
\author[1]{Romain Caron}
\author[1]{Fabrice Silva}
\author[1]{Christophe Vergez}
\author[1]{Bruno Cochelin}
\affil[1]{Aix Marseille Univ, CNRS, Centrale Med, LMA, Marseille, France}

\date{}

\begin{document}

\maketitle

\abstract{We present a minimal control-based continuation algorithm designed to track branches of limit cycles in autonomous systems. The controller can be viewed as three sub-controllers: (i) a derivative feedback controller that is used to stabilize the limit cycle, (ii) an integral phase controller, used to synthesize the unknown phase of the limit cycle and (iii) an integral arclength controller, used to track branches of limit cycles. The controlled system is analyzed theoretically, using the averaging method, allowing us to express tuning rules for the different parameters of the controller. Remarkably, theses tuning rules are independent of the studied system.}

\section{Introduction}

Bifurcation diagrams is a key tool to make a breach into the complex behavior of nonlinear system. They are used to explore how the behavior of a nonlinear dynamical system evolves when varying its parameters.  When a system model is available, numerical continuation methods are one of the standard tool to obtain these bifurcation diagrams. The fundamental aspect that allows to track branches of solutions of nonlinear dynamical system is the path following capability of numerical continuation algorithms. Different path following algorithm are available, such as pseudo-arclength parametrization, or the asymptotic numerical method, to name a few.
The idea of extending continuation method to nonlinear dynamics experiments has been first proposed by Sieber in \cite{sieber2007control}. All the necessary ingredients needed for experimental continuation are remarkably exposed in this seminal paper. They used a time delay feedback controller and introduced a dynamic phase condition to determine the period of the limit cycle, as well as a pseudo arclength condition to track these limit cycles. They applied numerically this control algorithm to track the limit cycles of a model of semiconductor laser with optical injection. In \cite{sieber2008control}, the same authors proposed to discretize the control signal by using Fourier series to increase the robustness of the method against potentially low accuracy of experimental measurements. The resulting fixed point problem is solved using quasi-Newton iteration and is applied numerically to the tracking of periodic orbit of a friction oscillator. A first experimental validation of control based continuation is presented in \cite{sieber2008experimental}, where it is used to track spinning orbits of a vertically excited pendulum. A similar control algorithm is used in \cite{barton2012control} and applied to the tracking of periodic solution of an energy harvester having both softening and hardening nonlinearities and in \cite{bureau2013experimental} to track periodic solution of a flexible beam with rigid constrains. 

Except methods using a time delay feedback controller, the above mentioned methods all rely on Newton iterations that require an estimate of the Jacobian matrix, which can be experimentally difficult due to the presence of noise. It was first proposed in \cite{barton2013systematic} to take the amplitude of the first Fourier coefficient as the continuation parameter. By doing so, making the controller non-invasive only requires to vanish the non-fundamental Fourier coefficient of the control output, which can be done by simple fixed point iteration, speeding up the process. Additionally, when applied to duffing-like systems, this method does not require any path following strategy since the response increases monotonously with respect to the forcing amplitude, which is also the continuation parameter. Since then, this method has been applied to the tracking of limit point of an energy harvester \cite{renson2017experimental}, to the tracking of periodic solution of a beam with harmonically coupled modes \cite{renson2019application} or to a system with friction nonlinearity \cite{kleyman2020application}.

As mentioned above, control based continuation methods have been applied successfully to various non-autonomous systems. Conversely, application to theses methods to autonomous systems has received less attention. An additional difficulty when dealing with non-autonomous systems is that the frequency of the limit cycle is also an unknown. Up to the author knowledge, the first experimental application of control based continuation to autonomous system is reported in \cite{tartaruga2019experimental}. The author used a control based continuation algorithm to obtain the bifurcation diagram of the limit cycles of an airfoil in a wind tunnel. To overcome the difficulty coming from the unknown frequency, the targeted limit cycle is parametrized by an angle in the phase space, instead of time. A similar approach, coupled to an adaptive B-spline parametrization is presented in \cite{blyth2023numerical}. The B-spline is used to provide a parsimonious representation of complex responses, allowing to track relaxation oscillations.

A point which still differentiates numerical continuation algorithm from experimental continuation ones, is the path following ability, that is, the capability of the controller to follow branches of solutions, even at turning point. Most of the time, this difficulty is overcome by using a natural continuation parameter, i.e. a continuation parameter that evolves monotonously. This is for example the case of phase locked loop based algorithm, where the complete bifurcation diagram of a Duffing oscillator can be obtained by imposing a monotonously increasing phase shift between the response and the excitation. The idea of adding an additional controller to impose an additional pseudo-arclength constrain has been first proposed by Siettos in \cite{siettos2004coarse}. The similar idea has been extended to periodic orbits in \cite{sieber2007control}. Recently, Raze proposed an arclength parametrization in \cite{raze2024experimental}. This method was successfully applied to various systems, including a plate exhibiting nonlinear modal interactions.

This paper aims at presenting a simple control based continuation algorithm able to track branches of limit cycles in autonomous systems. The designed controller can be viewed as an assembly of three already known controllers, namely: a derivative feedback controller, a phase locked loop controller and an arclength controller. The phase locked loop controller is used to synthesize the unknown phase of the limit cycle, and the arclength controller to allow path following capability. Even in their minimal form, control based continuation algorithms require the tuning of several control gains. The principal objective of this paper is to provide general tuning rules for this controller.

The paper is organized as follows. The designed controller is presented in Section \ref{sec:1}. The dynamics of the controlled system is analyzed theoretically using the method of averaging and is presented in Section \ref{sec:2}. The derivation of the stability conditions allows us to formulate generic tuning rules for the control gains. The behavior of the controller is numerically demonstrated on a generalized Van der Pol oscillator in Section \ref{sec:3}.

\section{Problem statement and design of the controller}\label{sec:1}

We consider a general weakly nonlinear single degree of freedom oscillator of the form

\begin{equation}\label{eq:1}
\ddot{x}+x=\epsilon g(x,\dot{x},\mu),
\end{equation}

where the dots represent time derivatives, $x$ is the displacement, $\mu$ is the bifurcation parameter, $\epsilon\ll1$ is a small parameter and $g(x,\dot{x},\mu)$ a smooth function. In this paper, we make the assumption that $g(x,\dot{x},\mu)$ is an odd function.

The objective is to design a controller that would allow us to track limit cycle oscillations using only the measured displacement $x$. This controller must satisfy the following requirements:

\begin{itemize}
\item Automatically converge to a limit cycle of the system
\item Stabilize the limit cycle if it is unstable
\item Non-invasivity, i.e. a limit cycle of the controlled system must be a limit cycle of the uncontrolled system
\item Path following capability.
\end{itemize}

The proposed controller can be splited into three sub-controllers described bellow.

\subsection{Derivative feedback controller}

According to control based continuation strategies, we consider a derivative feedback controller, such that Eq. (\ref{eq:1}) now reads

\begin{equation}\label{eq:2}
\ddot{x}+x=\epsilon g(x,\dot{x},\mu)+K_{1d}\dot{e}_1(t)
\end{equation}

where $K_{1d}$ is the derivative feedback gain and the error signal $e_1(t)$ is expressed by

\begin{equation}\label{eq:3}
e_1(t)=u(t)-x(t)
\end{equation}

Where $u(t)$ is the yet unknown control target. Clearly, if $u(t)=x(t)$, the control action is zero and the controller is non-invasive. Since we are considering weakly nonlinear system, we will assume that the response of the system is largely dominated by its fundamental harmonic such that we consider a single harmonic control target signal of the form

\begin{equation}\label{eq:4}
u(t)=G(t)\sin\theta(t)
\end{equation}

where $G(t)$ and $\theta(t)$ are the unknown amplitude and phase of the targeted limit cycle.

\subsection{Phase-locked loop controller}

An additional difficulty of the application of control-based continuation methods to autonomous system is that the period of limit cycle is also an unknown. This difficulty is addressed by using a phase-locked loop (PLL) controller that is used to synthesize the phase of the control signal. The PLL controller is composed of two parts: (i) a phase detector, (ii) an integral controller.

The phase detection is achieved by using synchronous demodulation as follow

\begin{equation}\label{eq:5}
\begin{array}{l}
\dot{y}_1=\omega_c\left(x\cos\theta-y_1\right)\\
\dot{y}_2=\omega_c\left(x\sin\theta-y_2\right)
\end{array}
\end{equation} 

where $y_1$ and $y_2$ are the output of the phase detector and $\omega_c$ is the cutoff frequency of the low pass filter. To understand the behavior of the phase detector, let us consider that the input displacement can be expressed by $x(t)=a\sin(\theta+\alpha)$. The two inputs signals of the phase detectors $x\cos\theta$ and $x\sin\theta$ are therefore expressed by

\begin{equation}\label{eq:6}
\begin{array}{l}
x\cos\theta=\dfrac{a}{2}\left(sin\alpha+\sin\left(2\theta+\alpha\right)\right)\vspace{0.1cm}\\
x\sin\theta=\dfrac{a}{2}\left(cos\alpha-\cos\left(2\theta+\alpha\right)\right)
\end{array}
\end{equation}

Assuming that the phase shift $\alpha$ evolves slowly compared to the instantaneous phase $\theta$ and perfect low pass filter, the output of the phase detector reads

\begin{equation}\label{eq:7}
y_1=\frac{a}{2}\sin\alpha,\quad y_2=\frac{a}{2}\cos\alpha
\end{equation}

Therefore, since the objective is to synthesize the control target signal such that it is in phase with the displacement of the oscillator, i.e. $\alpha=0$, the phase invasiveness error simply reads

\begin{equation}\label{eq:8}
e_2(t)=y_1(t)
\end{equation}

An integral controller is used to mitigate this error

\begin{equation}\label{eq:9}
\begin{array}{l}
\dot{y_3}=Re_2\\
\dot{\theta}=\omega_0+K_{i2}y_3
\end{array}
\end{equation}

where $\omega_0$ is the initial frequency and $K_{i2}$ is the integral control gain. The parameter $R$ could be omitted and is only used for scaling purpose.

\subsection{Arc controller}

The last building block of the controller aims at determining the amplitude of the control target $G(t)$ as well as allowing branch tracking capability. As seen in Eq. (\ref{eq:7}) if the PLL controller is successful (i.e. $\alpha=0$), $y_2=a/2$, such that the second output of the phase detector $y_2$ can be used as an estimate of the amplitude of the limit cycle. We can therefore define the amplitude invasiveness error as

\begin{equation}\label{eq:10}
e_3(t)=2y_2(t)-G(t)
\end{equation}

In order to allow branch tracking capability, similarly to numerical continuation, the amplitude $G(t)$ and the bifurcation parameter $\mu$ must be adjusted simultaneously. Similarly to \cite{raze2024experimental}, both the amplitude of the target signal $G(t)$ and the bifurcation parameter $\mu(t)$ are parametrized by and angle $\eta(t)$ on a circle of radius $\Delta$ centered at $(\mu_0,G_0)$ as 

\begin{equation}\label{eq:11}
\mu(t)=\mu_0+\Delta\cos\eta(t),\quad G(t)=G_0+\Delta\sin\eta(t),
\end{equation}

Notice that $\mu_0$, $G_0$ and $\Delta$ are additional user supplied parameters. Substituting Eq. (\ref{eq:11}) into Eq. (\ref{eq:10}) gives

\begin{equation}\label{eq:12}
e_3(t)=2y_2(t)-G_0-\Delta\sin\eta(t)
\end{equation}

A simple integral controller with gain $K_{i3}$ is used to mitigate the error $e_3(t)$ as

\begin{equation}\label{eq:13}
\dot{\eta}=K_{i3}e_3
\end{equation}

\subsection{Summary}

Gathering previous equations (\ref{eq:2}-\ref{eq:5},\ref{eq:9},\ref{eq:12},\ref{eq:13}), the equation of motion of the controlled nonlinear oscillator are expressed by

\begin{equation}\label{eq:14}
\begin{array}{l}
\ddot{x}+x=\epsilon g(x,\dot{x},\mu(\eta))+K_{d1}\left(\dot{\theta}\cos\theta\left(G_0+\Delta\sin\eta\right)+\Delta\dot{\eta}\cos\eta\sin\theta-\dot{x}\right)\\
\dot{y}_1=\omega_c\left(x\cos\theta-y_1\right)\\
\dot{y}_2=\omega_c\left(x\sin\theta-y_2\right)\\
\dot{y}_3=Ry_1\\
\dot{\theta}=\omega_0+K_{i2}y_3\\
\dot{\eta}=K_{i3}\left(2y_2-G_0-\Delta\sin\eta\right)
\end{array}
\end{equation}

\section{Theoretical analysis}\label{sec:2}

\subsection{Averaged equations}

In this section, the objective is to express the slow flow equations of Eq. (\ref{eq:14}) by using the method of averaging. The displacement and velocity of the oscillator are expressed by

\begin{equation}\label{eq:15}
x(t)=a(t)\sin\phi(t),\quad \dot{x}(t)=a(t)\cos\phi(t)
\end{equation}

where $a(t)$ is the amplitude of the oscillations and $\phi(t)=\theta(t)+\alpha(t)$ is the instantaneous phase. Expressing time derivatives of Eq. (\ref{eq:15}) 

\begin{equation}\label{eq:16}
\begin{array}{l}
\dot{x}=\dot{a}\sin\phi+a\dot{\phi}\cos\phi\\
\ddot{x}=\dot{a}\cos\phi-a\dot{\phi}\sin\phi
\end{array}
\end{equation}

Equating the second equation of Eq. (\ref{eq:15}) and the first equation of Eq. (\ref{eq:16}) gives

\begin{equation}\label{eq:17}
a\left(1-\dot{\theta}-\dot{\alpha}\right)\cos\phi-\dot{a}\sin\phi=0
\end{equation}

The method of averaging apply if the time derivatives evolves slowly, i.e. $d/dt\sim\mathcal{O}(\epsilon)$. Accordingly, the control parameters are rescaled with respect to the small parameter $\epsilon$ to explicit the condition on the orders of magnitude as follow

\begin{equation}\label{eq:18}
\left\lbrace\omega_c,R,K_{d1},K_{i2},K_{i3}\right\rbrace=\epsilon\left\lbrace\bar{\omega}_c,\bar{R},\bar{K}_{d1},\bar{K}_{i2},\bar{K}_{i3}\right\rbrace
\end{equation}

In addition, we consider that the initial frequency of the PLL controller $\omega_0$ is close to the frequency of the emerging limit cycle, i.e. $\omega_0=1+\epsilon\sigma$.

Equations (\ref{eq:15}, \ref{eq:16}) are substituted into the extended system Eq. (\ref{eq:14},\ref{eq:17}). Solving for first order derivative terms and keeping only terms up to $\mathcal{O}(\epsilon)$ gives

\begin{equation}\label{eq:19}
\begin{array}{l}
\dot{a}=\epsilon\left[g(a\sin\phi,a\cos\phi,\mu(\eta))\cos\phi+\bar{K}_{d1}\cos\phi\left((G_0+\Delta\sin\eta)\cos(\phi-\alpha)-a\cos\phi\right)\right]\\
\begin{aligned}
\dot{\alpha}=-\dfrac{\epsilon}{a}\left[g(a\sin\phi,a\cos\phi,\mu(\eta))\sin\phi +a(\sigma+\bar{K}_{i2}y_3)+\right.\\
\left. \bar{K}_{d1}\sin\phi\left((G_0+\Delta\sin\eta)\cos(\phi-\alpha)-a\cos\phi\right)\right]
\end{aligned}\\
\dot{y}_1=\epsilon\bar{\omega}_c\left(a\sin\phi\cos(\phi-\alpha)-y_1\right)\\
\dot{y}_2=\epsilon\bar{\omega}_c\left(a\sin\phi\sin(\phi-\alpha)-y_2\right)\\
\dot{y}_3=\epsilon \bar{R}y_1\\
\dot{\eta}=\epsilon \bar{K}_{i3}\left(2y_2-G_0-\Delta\sin\eta\right)\\
\dot{\theta}=1+\epsilon\left(\sigma+\bar{K}_{i2}y_3\right)
\end{array}
\end{equation}

Averaged equations are obtained by integrating Eq. (\ref{eq:19}) over a period ($\frac{1}{2\pi}\int_0^{2\pi}d\phi$). Reabsorbing $\epsilon$ by making inverse scaling from Eq. (\ref{eq:18}) gives the following slow flow system

\begin{equation}\label{eq:20}
\begin{array}{l}
\dot{a}=\epsilon f_1(a,\mu(\eta))+\dfrac{K_{d1}}{2}\left((G_0+\Delta\sin\eta)\cos\alpha-a\right)\\
\dot\alpha=-\dfrac{\epsilon}{a}f_2(a,\mu(\eta))+1-\omega_0-K_{i2}y_3-\dfrac{K_{d1}}{2a}(G_0+\Delta\sin\eta)\sin\alpha\\
\dot{y}_1=\omega_c\left(\dfrac{a}{2}\sin\alpha-y_1\right)\\
\dot{y}_2=\omega_c\left(\dfrac{a}{2}\cos\alpha-y_2\right)\\
\dot{y}_3=Ry_1\\
\dot{\eta}=K_{i3}\left(2y_2-G_0-\Delta\sin\eta\right)\\
\dot{\theta}=\omega_0+K_{i2}y_3
\end{array}
\end{equation}

The only system-dependent functions $f_1$ and $f_2$ are expressed by

\begin{equation}\label{eq:21}
\begin{array}{l}
f_1(a,\mu(\eta))=\dfrac{1}{2\pi}\int_{0}^{2\pi} g(a\sin\phi,a\cos\phi,\mu(\eta))\cos\phi d\phi\\
f_2(a,\mu(\eta))=\dfrac{1}{2\pi}\int_{0}^{2\pi} g(a\sin\phi,a\cos\phi,\mu(\eta))\sin\phi d\phi
\end{array}
\end{equation}

Notice that since $g$ is an odd function, $f_1$ and $f_2$ does not depend on $\alpha$.

\subsection{Uncontrolled system}

Before analyzing the behavior of the controlled system it is useful to describe the behavior of the uncontrolled system. The averaged equations of the uncontrolled, obtained by setting the control gains to zero in Eq. (\ref{eq:20}), are given by

\begin{equation}\label{eq:22}
\begin{array}{l}
\dot{a}=\epsilon f_1(a,\mu)\\
\dot\alpha=-\dfrac{\epsilon}{a}f_2(a,\mu)+1-\omega_0\\
\dot{\theta}=\omega_0
\end{array}
\end{equation}

Denoting $\Omega=\dot{\theta}$, the instantaneous frequency of the limit cycle, the fixed points $(\tilde{a},\tilde{\alpha})$ of Eq. (\ref{eq:22}) are obtained by setting the corresponding time derivatives to zero yielding

\begin{equation}\label{eq:23}
\begin{array}{l}
f_1(\tilde{a},\mu)=0\\
\omega_0=1-\dfrac{\epsilon}{\tilde{a}}f_2(\tilde{a},\mu)\\
\dot{\theta}\equiv\Omega=\omega_0
\end{array}
\end{equation}

Notice that the phase $\tilde{\alpha}$ is undetermined, consistently with the fact that the system is autonomous. Combining the last two equations in Eq. (\ref{eq:23}), the instantaneous frequency of the limit cycle of the uncontrolled system reads

\begin{equation}\label{eq:24}
\Omega=1-\frac{\epsilon}{\tilde{a}}f_2(\tilde{a},\mu)
\end{equation}

The stability of the fixed point is obtained by linearizing the system around the fixed point and is given by

\begin{equation}\label{eq:25}
\lambda_u=\epsilon f_{1a}(\tilde{a},\mu)
\end{equation}

where $f_{jq}=\partial f_j/\partial q$. The limit cycle is unstable if $\lambda_u>0$.

\subsection{Controlled system}

\subsubsection{Fixed points}

The fixed points of the controlled system $(\tilde{a},\tilde{\alpha},\tilde{y}_1,\tilde{y}_2,\tilde{y}_3,\tilde{\eta})$ are obtained by equating the right hand side of Eq. (\ref{eq:20}) to zero, giving

\begin{equation}\label{eq:26}
\begin{array}{l}
0=\epsilon f_1(\tilde{a},\mu(\tilde{\eta}))+\dfrac{K_{d1}}{2}\left((G_0+\Delta\sin\tilde{\eta})\cos\tilde{\alpha}-\tilde{a}\right)\\
0=-\dfrac{\epsilon}{\tilde{a}}f_2(\tilde{a},\mu(\tilde{\eta}))+1-\omega_0-K_{i2}\tilde{y}_3-\dfrac{K_{d1}}{2\tilde{a}}(G_0+\Delta\sin\tilde{\eta})\sin\tilde{\alpha}\\
0=\omega_c\left(\dfrac{\tilde{a}}{2}\sin\tilde{\alpha}-\tilde{y}_1\right)\\
0=\omega_c\left(\dfrac{\tilde{a}}{2}\cos\tilde{\alpha}-\tilde{y}_2\right)\\
0=R\tilde{y}_1\\
0=K_{i3}\left(2\tilde{y}_2-G_0-\Delta\sin\tilde{\eta}\right)\\
\Omega=\omega_0+K_{i2}\tilde{y}_3
\end{array}
\end{equation}

As is Eq. (\ref{eq:23}), we have introduced the instantaneous frequency $\Omega=\dot{\theta}$. Substituting the expression from arclength parametrization (Eq. (\ref{eq:11})) $\sin\eta=(G(\eta)-G_0)/\Delta$ into Eq. (\ref{eq:26}) and solving gives

\begin{equation}\label{eq:27}
\begin{array}{l}
\sin\tilde{\alpha}=0\\
\cos\tilde{\alpha}=\dfrac{G}{\tilde{a}}\\
\tilde{y}_1=0\\
\tilde{y}_2=\dfrac{\tilde{a}}{2}\\
f_1(\tilde{a},\mu(\tilde{\eta}))=\dfrac{K_{d1}\pi}{\tilde{a}\epsilon}(\tilde{a}^2-G(\tilde{\eta})^2)\\
\tilde{y}_3=\dfrac{1}{K_{i2}}\left(1-\omega_0-\dfrac{\epsilon}{\tilde{a}}f_2(\tilde{a},\mu(\tilde{\eta}))\right)\\
\dot{\theta}\equiv\Omega=1-\dfrac{\epsilon}{\tilde{a}}f_2(\tilde{a},\mu(\tilde{\eta}))
\end{array}
\end{equation}

Considering only positive amplitudes, i.e. $G>0, a>0$ the fixed points are finally expressed as

\begin{equation}\label{eq:28}
\begin{array}{l}
f_1(\tilde{a},\mu(\tilde{\eta}))=0\\
\tilde{\alpha}=0\\
\Omega=1-\dfrac{\epsilon}{\tilde{a}}f_2(\tilde{a},\mu(\tilde{\eta}))\\
\cos\tilde{\eta}=\dfrac{\mu-\mu_0}{\Delta},\quad \sin\tilde{\eta}=\dfrac{\tilde{a}-G_0}{\Delta}\\
G(\tilde{\eta})=\tilde{a}\\
\tilde{y}_1=0,\quad \tilde{y}_2=\dfrac{\tilde{a}}{2},\quad \tilde{y}_3=\dfrac{\Omega-\omega_0}{K_{i2}}
\end{array}
\end{equation}

Therefore, if they exist, the fixed points of the controlled system corresponds to fixed points of the uncontrolled system, located at a distance $\Delta$ from $(\mu_0,G_0)$, i.e. the controller is non-invasive.

\subsubsection{Stability}

The stability of the fixed points is determined by looking at the eigenvalues of the following characteristic matrix obtained by linearizing Eq. (\ref{eq:20}) around the fixed point

\begin{equation}\label{eq:29}
\bm{M}=\begin{pmatrix}
-\dfrac{K_{d1}}{2}+\epsilon f_{1a} & 0 & 0 & 0 & 0 & \dfrac{K_{d1}}{2}(\tilde{\mu}-\mu_0)+\epsilon f_{1\mu}\dfrac{d\mu}{d\eta}\\
\dfrac{\epsilon}{\tilde{a}^2}(\tilde{a}f_{2a}-f_2) & -\dfrac{K_{d1}}{2} & 0 & 0 & -K_{i2} & -\dfrac{\epsilon}{\tilde{a}}f_{2\mu}\dfrac{d\mu}{d\eta}\\
0 & \dfrac{\tilde{a}\omega_c}{2} & -\omega_c & 0 & 0 & 0\\
\dfrac{\omega_c}{2} & 0 & 0 & -\omega_c & 0 & 0\\
0 & 0 & R & 0 & 0 & 0\\
0 & 0 & 0 & 2K_{i3} & 0 & -K_{p3}(\tilde{\mu}-\mu_0)
\end{pmatrix}
\end{equation}

where $\tilde{\mu}=\mu(\tilde{\eta})$. The corresponding characteristic polynomial can be factorized into two cubic polynomials in $\lambda$ expressed by

\begin{equation}\label{eq:30}
\begin{array}{l}
0=2\lambda^3+(K_{d1}+2\omega_c)\lambda^2+\omega_c K_{d1}\lambda+\omega_c\tilde{a}RK_{i2}\\
\begin{aligned}
0=2\lambda^3+(2\omega_c+K_{d1}+2K_{i3}(\tilde{\mu}-\mu_0)-2\epsilon f_{1a})\lambda^2\\
+(\omega_c(K_{d1}-2\epsilon f_{1a})+K_{i3}(\tilde{\mu}-\mu_0)(2\omega_c+K_{d1}-2\epsilon f_{1a}))\lambda\\
-2K_{i3}\epsilon\omega_c\left((\tilde{\mu}-\mu_0)f_{1a}+\dfrac{d\mu}{d\eta}f_{1\mu}\right)
\end{aligned}
\end{array}
\end{equation}

We first deal with the first equation of Eq. (\ref{eq:30}). Remarkably, it is not system dependent since neither $f_1$ or $f_2$ appears in this equation. The stability condition can be obtained by using Routh-Hurwitz criterion as follow

\begin{equation}\label{eq:31}
\begin{array}{l}
\dfrac{K_{d1}}{2}+\omega_c>0\\
K_{d1}>0\\
RK_{i2}>0\\
(K_{d1}+2\omega_c)K_{d1}-\tilde{a}RK_{i2}>0
\end{array}
\end{equation}

The first three equations are trivially satisfied if the gains are positive, i.e. $K_{d1}>0$, $K_{i2}>0$ and $R>0$. Isolating $\tilde{a}$ from the fourth equation yields to

\begin{equation}\label{eq:32}
\tilde{a}<\frac{K_{d1}(K_{d1}+2\omega_c)}{2RK_{i2}}
\end{equation}

It follows from Eq. (\ref{eq:32}) that the controller will become unstable if this condition is violated. In other words, we are able to determine, independently of the studied system, the maximal amplitude of the limit cycle that we will be able to stabilize.

Let us now consider the second characteristic polynomial in Eq. (\ref{eq:30}). Considering the expression of the fixed points, the derivative of $\mu$ with respect to $\eta$ is expressed as follow

\begin{equation}
\left.\dfrac{d\mu}{d\eta}\right|_{\tilde{\mu},\tilde{\eta}}=-(\tilde{a}-G_0)
\end{equation}

Substituting into Eq. (\ref{eq:30}) and considering small $\Delta$ such that $(a-G_0)\ll1$ and $(\mu-\mu_0)\ll1$ a perturbative solution can be expressed by

\begin{equation}\label{eq:34}
\begin{array}{l}
\lambda_1=-\dfrac{2\epsilon K_{i3}}{2\epsilon f_{1a}-K_{d1}}\left((\tilde{\mu}-\mu_0)f_{1a}-(a-G_0)f_{1\mu}\right)+\mathcal{O}(\Delta^2)\\
\lambda_2=-\omega_c+\mathcal{O}(\Delta)\\
\lambda_3=\epsilon f_{1a}-\dfrac{K_{d1}}{2}+\mathcal{O}(\Delta)
\end{array}
\end{equation}

Remembering Eq. (\ref{eq:25}) (i.e. $\epsilon f_{1a}=\lambda_u$ is the eigenvalue of the uncontrolled system), $\lambda_3$ in Eq. (\ref{eq:34}) shows that unstable limit cycle of the uncontrolled system can be stabilized by choosing sufficiently large $K_{d1}$. Le us now analyze $\lambda_1$ in Eq. (\ref{eq:34}). The normal vector to the bifurcation diagram at the fixed point is $\bm{n}=(f_{1\mu},f_{1a})^T$. The tangent vector $\bm{t}$, satisfying $\bm{n}^T\bm{t}=0$ can be expressed by

\begin{equation}\label{eq:35}
\bm{t}=\begin{pmatrix}
-f_{1a}\\f_{1\mu}
\end{pmatrix}
\end{equation}

Therefore, $\lambda_1$ in Eq. (\ref{eq:34}) can be rewritten as follow

\begin{equation}\label{eq:36}
\lambda_1=\dfrac{2\epsilon K_{i3}}{2\epsilon f_{1a}-K_{d1}}\bm{t}^T\begin{pmatrix}
\tilde{\mu}-\mu_0\\ a-G_0
\end{pmatrix}
\end{equation}

\begin{figure}
\begin{center}
\includegraphics[scale=1]{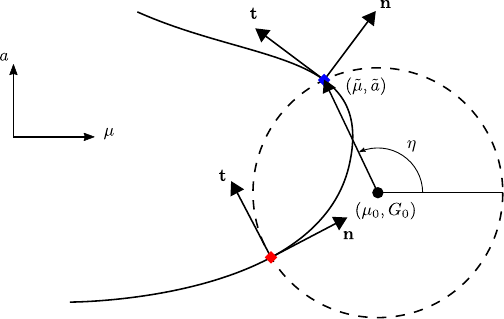}
\caption{Graphical interpretation of the stability condition given in Eq. (\ref{eq:36}). Blue and red squares correspond to stable and unstable fixed point, respectively.}
\label{fig:1}
\end{center}
\end{figure}

Assuming that the derivative control gain $K_{d1}$ is such that $\lambda_3<0$ (i.e. the limit cycle can be stabilized), it follows that the denominator in Eq. (\ref{eq:36}) is negative. Therefore, if the scalar product in Eq. (\ref{eq:36}) is positive, the fixed point will be stable if the proportional gain $K_{i3}>0$. On the contrary, if the scalar product in Eq. (\ref{eq:36}) is negative, the fixed point will be stable for $K_{i3}<0$. A graphical interpretation is depicted in Fig. \ref{fig:1}. Assuming that $K_{i3}>0$ the stable fixed point, corresponding to positive scalar product is the blue one, while the red one is unstable.

Let us consider that the $(\mu_0,G_0)$ correspond to an already converged limit cycle, i.e. $(\mu_0,G_0)$ is a fixed point of the uncontrolled system. If the radius $\Delta$ is small enough, the circle intersects twice the bifurcation diagram such that there exist two fixed point, then
\begin{itemize}
\item The direction of continuation is chosen by the sign of $K_{p3}$\\
\item The controller can handle turning points.
\end{itemize}

\section{Numerical example}\label{sec:3}

As an application example, we consider a generalized Van der Pol oscillator. The nonlinear function $g(x,\dot{x},\mu)$ is given by

\begin{equation}\label{eq:37}
g(x,\dot{x},\mu)=\left(\mu+\beta x^2-x^4\right)\dot{x}-\rho x^3
\end{equation}

Using Eq. (\ref{eq:21}), $f_1$ and $f_2$ are expressed by

\begin{equation}\label{eq:38}
\begin{array}{l}
f_1(a,\mu)=\dfrac{a}{16}(\mu+2a^2\beta-a^2)\\
f_2(a,\mu)=-\dfrac{3}{8}\rho a^3
\end{array}
\end{equation}

The fixed points of the system are obtained using Eq. (\ref{eq:23}) and are expressed as follow

\begin{equation}\label{eq:39}
\tilde{a}=\sqrt{\beta\pm\sqrt{\beta^2+8\tilde{\mu}}}
\end{equation}

\begin{figure}
\begin{center}
\includegraphics[scale=1]{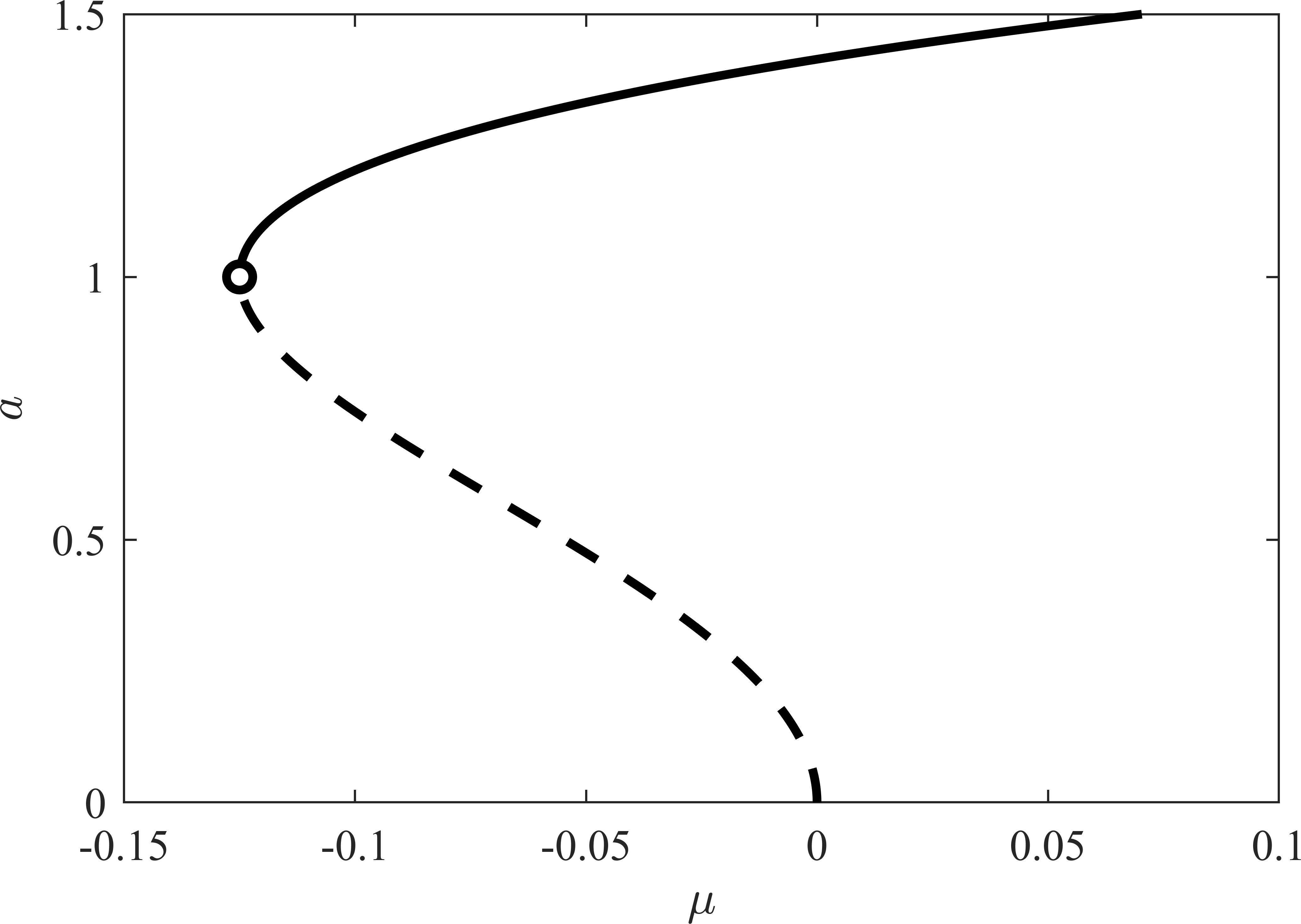}
\caption{Bifurcation diagram of the uncontrolled system Eq. (\ref{eq:37}) with $\beta=1$. plain and dotted lines represent the amplitude of stable and unstable limit cycle, respectively. The dot corresponds to fold bifurcation.}
\label{fig:2}
\end{center}
\end{figure}

A fold bifurcation occurs at $\tilde{\mu}=-\beta^2/8$. An example of bifurcation diagram for $\beta=1$ is depicted in Fig. \ref{fig:2}.

We will now investigate the behavior of the controlled system. The parameters of the controller, chosen according to the previously obtained tuning rules, are gathered in Table \ref{tab:1}. The starting point has been chosen arbitrarily as $(\mu_0,G_0)=(0,0.3)$. In order to highlight the raw behavior of the controller, all initial conditions are set to zero, which is obviously not an optimal choice for practical application.

\begin{table}
\begin{center}
\begin{tabular}{|l l|}
\hline
$K_{d1}=0.1$ & $\omega_c=0.01$\\
$K_{i2}=0.1$ & $\omega_0=0.9$\\
$K_{i3}=0.1$ & $\Delta=0.1$\\
$R=0.1$ & $\epsilon=0.1$\\
\hline
\end{tabular}
\caption{Parameters of the controller}
\label{tab:1}
\end{center}
\end{table}

\begin{figure}
\begin{center}
\includegraphics[scale=1]{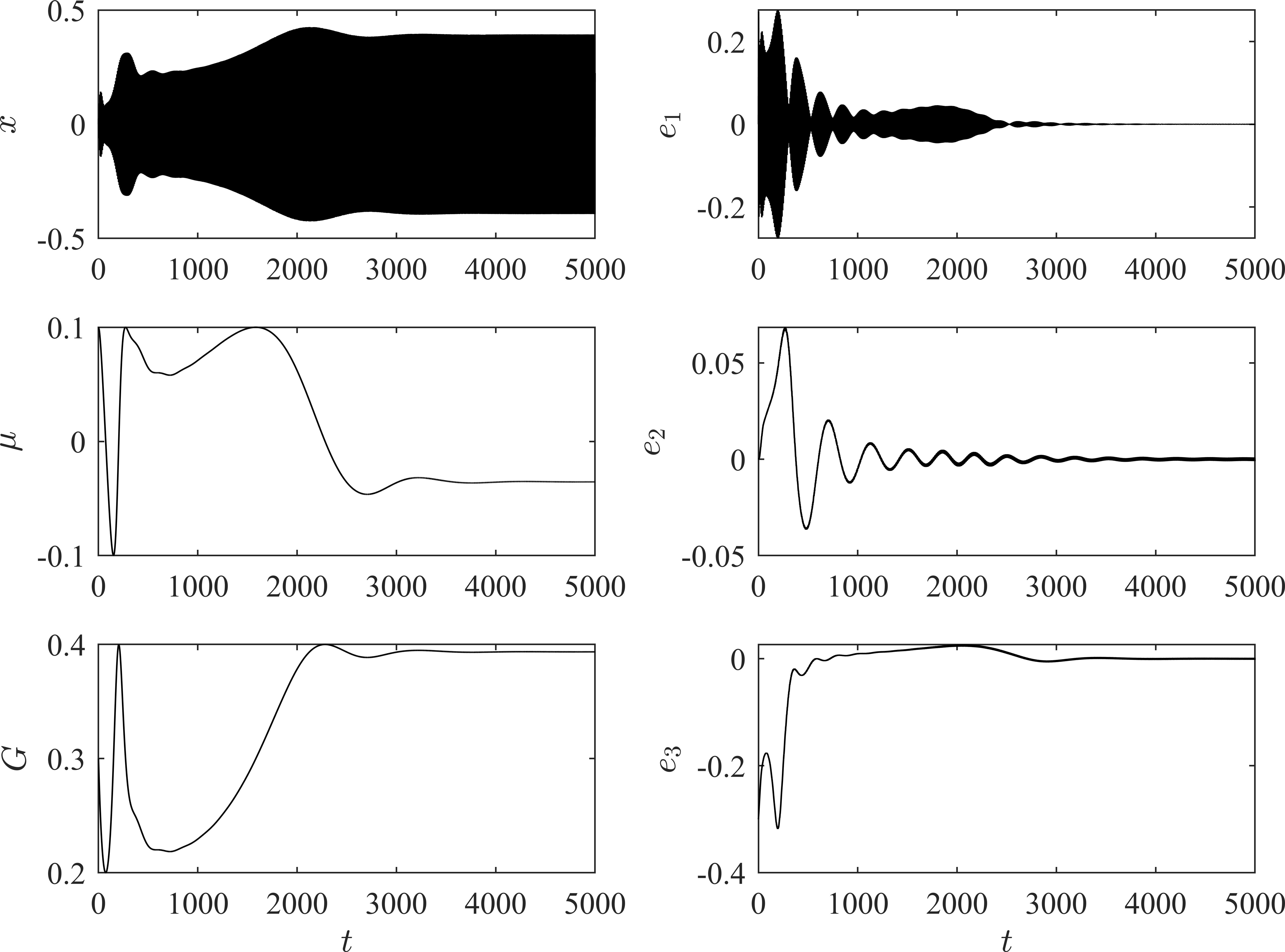}
\caption{Result of numerical integration for parameters in Table \ref{tab:1}.}
\label{fig:3}
\end{center}
\end{figure}

\begin{figure}
\begin{center}
\includegraphics[scale=1]{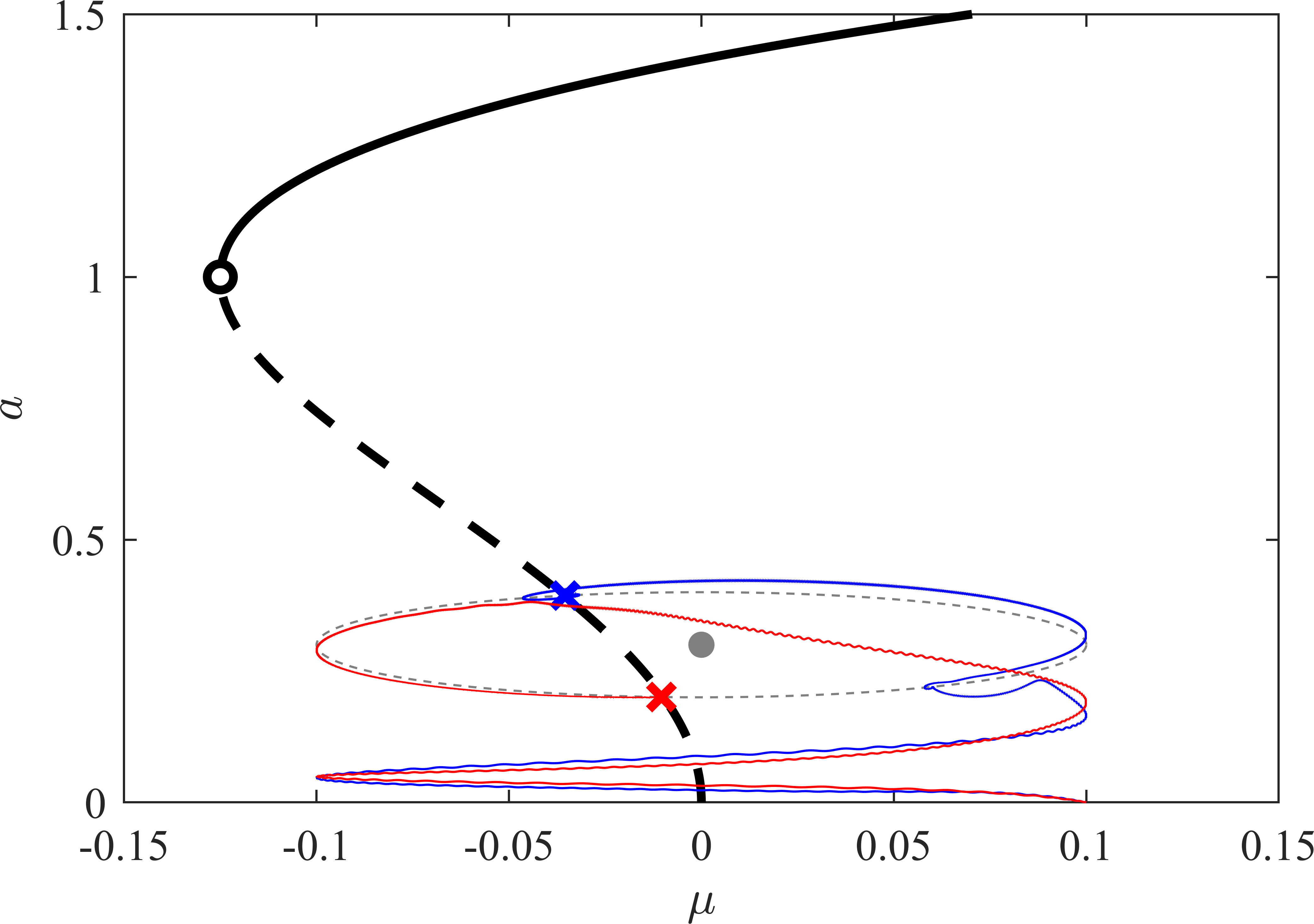}
\caption{Result of numerical integration for the bifurcation parameters $\mu$ and $G$ plotted on the bifurcation diagram of the uncontrolled system for $K_{p3}=0.1$ (blue) and $K_{p3}=-0.1$ (red). The other parameters are given in Table \ref{tab:1}. The grey dot represents the starting point $(\mu_0,G_0)=(0,0.3)$ and the dotted grey line, the circle of radius $\Delta=0.1$.} 
\label{fig:4}
\end{center}
\end{figure}

The result of numerical integration is depicted in Fig. \ref{fig:3}. It is clearly observed that the controller converges to a stable limit cycle, and we can confirm that this limit cycle is also a limit cycle of the uncontrolled system since the error signals converge to zero. The evolution of the bifurcation parameter $\mu$ and of the estimated amplitude of the limit cycle $\hat{a}=2y_2$ are also depicted on the bifurcation diagram in Fig. \ref{fig:4} and correspond to the blue line. The crosses are simply the last point of the numerical integration. The central point $(\mu_0,G_0)$ is depicted by the grey dot and the grey dashed line corresponds to the circle of diameter $\Delta$. It is observed that after some transient the flow converges to the unstable branch of the theoretical bifurcation diagram. Notice that the starting point is not necessarily a point on the bifurcation diagram, but has to be close enough such that the circle of diameter $\Delta$ intersects the bifurcation diagram.  

The same numerical integration has been performed for $K_{p3}=-0.1$ and the result in depicted by the red line in Fig. \ref{fig:4}. It this case, the stable fixed point is the lower one, in accordance with the described theoretical behavior.

\begin{figure}
\begin{center}
\includegraphics[scale=1]{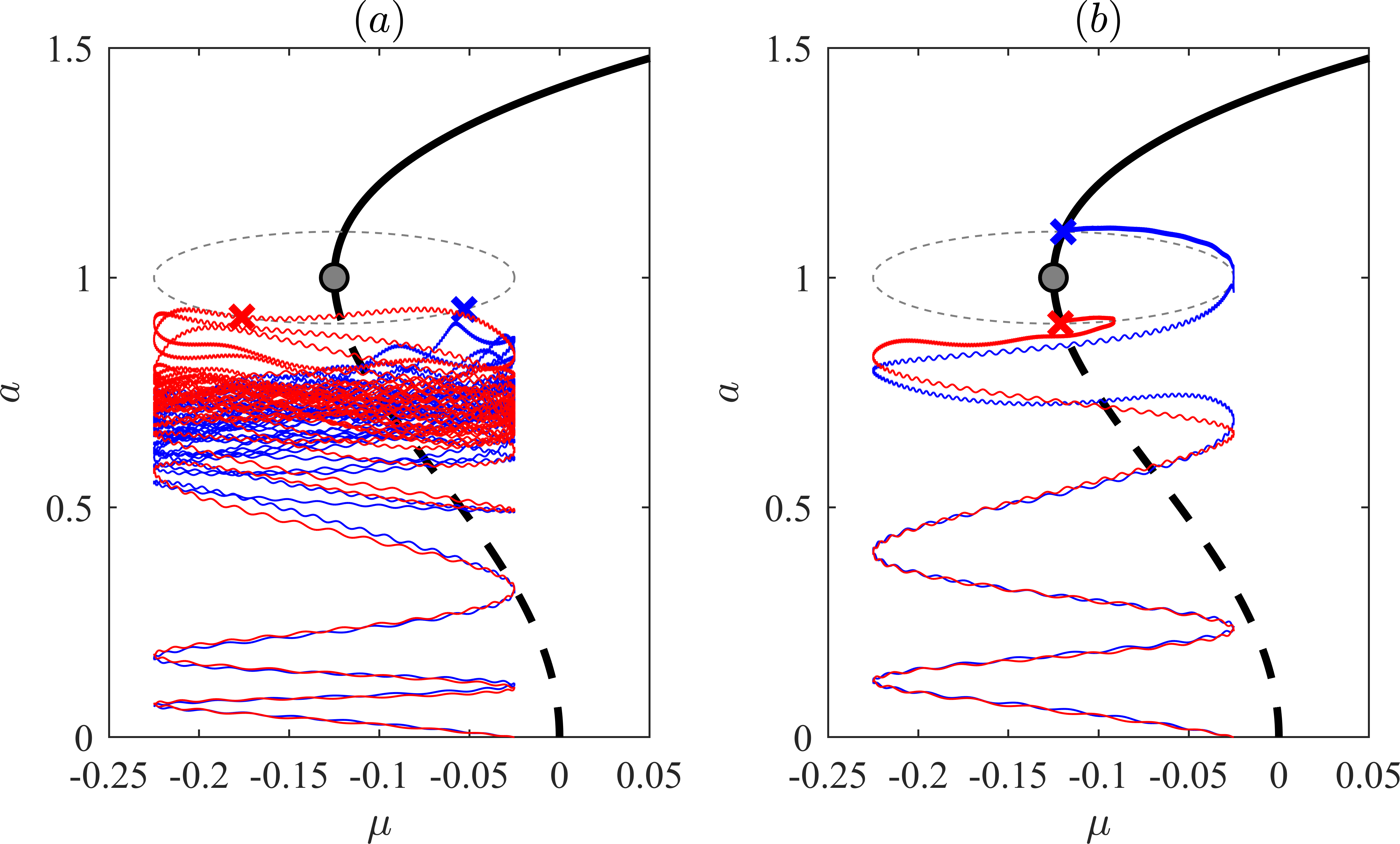}
\caption{Result of numerical integration of the controlled system displayed on the bifurcation diagram of the uncontrolled system for $K_{p3}=0.1$ (blue) and $K_{p3}=-0.1$ (red). (a) $K_{d1}=0.1$, (b) $K_{d1}=0.2$, the other parameters are given in Table \ref{tab:1}. The grey dot represents the starting point on the fold bifurcation point.}
\label{fig:5}
\end{center}
\end{figure}

A similar numerical experiment can be performed with a starting point located at the fold bifurcation.  Again, the blue line corresponds to the numerical integration with $K_{p3}=0.1$ while the red line corresponds to $K_{p3}=-0.1$. Notice that for the parameters given in Tab. \ref{tab:1}, the maximal amplitude of the limit cycle that the controller is able to stabilize, given by Eq. (\ref{eq:32}) is $a_{max}=0.6$. It follows that for the parameters given in Tab. \ref{tab:1}, the controller is not able to converge to the limit cycles close to the Fold bifurcation. This is confirmed numerically as depicted in Fig. \ref{fig:5}(a) where it is clearly observed that independently on the sign of $K_{i3}$, the system does not converge to the limit cycle. Increasing the derivative gain to $K_{d1}=0.2$, now $a_{max}=2.2$. In this case, and according to the tuning rule, the system converge to the limit cycle either above or under the Fold bifurcation point, depending on the sign of $K_{i3}$, as depicted in Fig.  \ref{fig:5}(b).

This clearly shows that thanks to the arclength parametrization, the control algorithm is able to handle the turning point. It also shows that the continuation direction is preserved even if the open loop stability of the periodic solution changes.

\begin{figure}
\begin{center}
\includegraphics[scale=1]{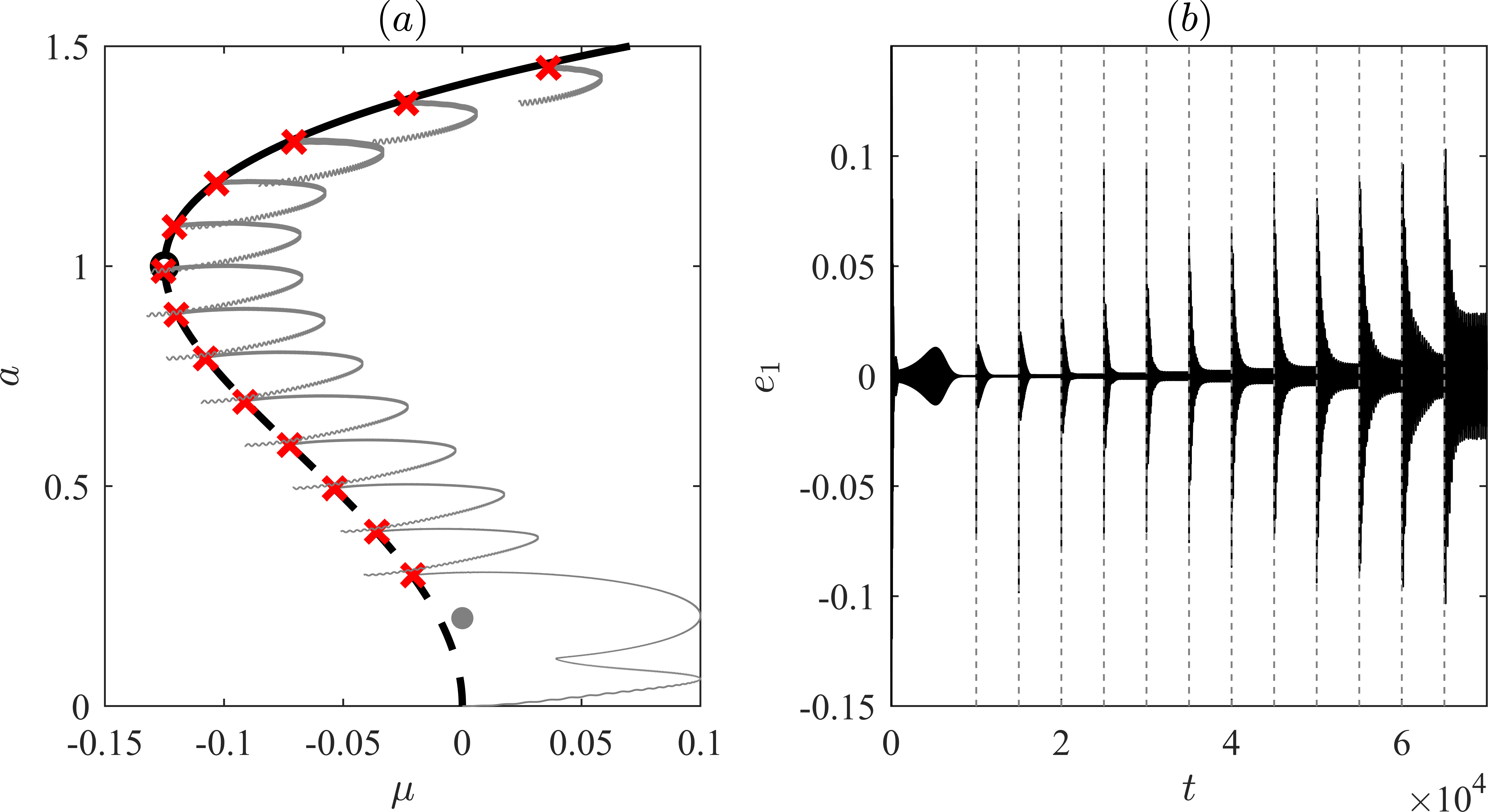}
\caption{(a): Reconstructed bifurcation diagram for $K_{i3}=0.1$ and $K_{d1}=0.2$. Red crosses indicates the value of the bifurcation parameter and the estimated amplitude of the limit cycle at the end of each numerical integration. (b): evolution of the error $e_1(t)$.}
\label{fig:6}
\end{center}
\end{figure}

The bifurcation diagram of the studied system can be reconstructed by successively waiting for transient die out and updating the central point $(\mu_0,G_0)$. The obtained bifurcation diagram, for $K_{i3}=0.1$ and $K_{d1}=0.2$ is depicted in Fig. \ref{fig:6}(a). The red crosses indicate the value of the bifurcation parameter and the estimated amplitude of the limit cycle at the end of each numerical integration. The evolution of the error $e_1(t)$ (given in Eq. (\ref{eq:3})), showing the invasiveness of the controller is depicted in Fig. \ref{fig:6}(b). It is clearly observed that the error increases when the amplitude of the limit cycle increases, which is attributed to the presence of higher harmonics that are not balanced by the mono-harmonic control target (Eq. (\ref{eq:4})), but remain small compared to the amplitude of the limit cycle.

\section{Conclusion}

We have presented a minimal control based continuation controller that is able to stabilize and follow branches of limit cycles. The controller can be described as an assembly of three sub-controllers:
\begin{itemize}
\item A derivative feedback controller aiming to stabilize the limit cycle.
\item An integral phase controller used to synthesize the unknown phase of the limit cycle.
\item An integral arclength controller used to track the limit cycle.
\end{itemize}

In order to derive tuning rules for this controller, the dynamics of the controlled system has been analyzed theoretically using the method of averaging. The analysis of the stability of the reduced slow flow system allowed us to derive the following tuning rules
\begin{itemize}
\item The derivative gain $K_{d1}$ and the integral gain $K_{i2}$ must be positive.
\item An unstable limit cycle can be stabilized by using a sufficiently large value of $K_{d1}$.
\item With given control gains and frequencies of the low pass filter, the maximal amplitude of the limit cycle that can be stabilized is given by Eq. (\ref{eq:32}).
\item The direction of continuation is chosen by the sign of the integral gain $K_{i3}$.
\end{itemize}

Remarkably, the tuning rules are independent of the studied system.
The designed control algorithm has been validated numerically on a generalized Van der Pol oscillator whose bifurcation diagram exhibit a cyclic fold bifurcation. These numerical experiments also confirmed the ability of the controller to deal with turning points.

The main limitation of the presented control algorithm is the limitation to mono-harmonic control target preventing the non-invasiveness of the controller when the system exhibits higher harmonics. Adaptive filters seems to be a good candidate to alleviate this limitation.

\bibliographystyle{plain}
\bibliography{biblio}

\begin{thebibliography}{10}

\bibitem{barton2012control}
D.~Barton, B.~Mann, and S.~Burrow.
\newblock Control-based continuation for investigating nonlinear experiments.
\newblock {\em Journal of Vibration and Control}, 18(4):509--520, 2012.

\bibitem{barton2013systematic}
D.~Barton and J.~Sieber.
\newblock Systematic experimental exploration of bifurcations with noninvasive
  control.
\newblock {\em Physical Review E—Statistical, Nonlinear, and Soft Matter
  Physics}, 87(5):052916, 2013.

\bibitem{blyth2023numerical}
M.~Blyth, K.~Tsaneva-Atanasova, L.~Marucci, and L.~Renson.
\newblock Numerical methods for control-based continuation of relaxation
  oscillations.
\newblock {\em Nonlinear Dynamics}, 111(9):7975--7992, 2023.

\bibitem{bureau2013experimental}
E.~Bureau, F.~Schilder, I-F. Santos, J-J. Thomsen, and J.~Starke.
\newblock Experimental bifurcation analysis of an impact oscillator—tuning a
  non-invasive control scheme.
\newblock {\em Journal of Sound and Vibration}, 332(22):5883--5897, 2013.

\bibitem{kleyman2020application}
G.~Kleyman, M.~Paehr, and S.~Tatzko.
\newblock Application of control-based-continuation for characterization of
  dynamic systems with stiffness and friction nonlinearities.
\newblock {\em Mechanics Research Communications}, 106:103520, 2020.

\bibitem{raze2024experimental}
G.~Raze, G.~Abeloos, and G.~Kerschen.
\newblock Experimental continuation in nonlinear dynamics: recent advances and
  future challenges.
\newblock {\em Nonlinear Dynamics}, pages 1--49, 2024.

\bibitem{renson2017experimental}
L.~Renson, D.~Barton, and S.~Neild.
\newblock Experimental tracking of limit-point bifurcations and backbone curves
  using control-based continuation.
\newblock {\em International Journal of Bifurcation and Chaos}, 27(01):1730002,
  2017.

\bibitem{renson2019application}
L.~Renson, A.~Shaw, D.~Barton, and S.~Neild.
\newblock Application of control-based continuation to a nonlinear structure
  with harmonically coupled modes.
\newblock {\em Mechanical Systems and Signal Processing}, 120:449--464, 2019.

\bibitem{sieber2008experimental}
J.~Sieber, A.~Gonzalez-Buelga, S.~Neild, D.~Wagg, and B.~Krauskopf.
\newblock Experimental continuation of periodic orbits through a fold.
\newblock {\em Physical review letters}, 100(24):244101, 2008.

\bibitem{sieber2007control}
J.~Sieber and B.~Krauskopf.
\newblock Control-based continuation of periodic orbits with a time-delayed
  difference scheme.
\newblock {\em International Journal of Bifurcation and Chaos},
  17(08):2579--2593, 2007.

\bibitem{sieber2008control}
J.~Sieber and B.~Krauskopf.
\newblock Control based bifurcation analysis for experiments.
\newblock {\em Nonlinear Dynamics}, 51(3):365--377, 2008.

\bibitem{siettos2004coarse}
C.~Siettos, I.~Kevrekidis, and D.~Maroudas.
\newblock Coarse bifurcation diagrams via microscopic simulators: a
  state-feedback control-based approach.
\newblock {\em International Journal of Bifurcation and Chaos},
  14(01):207--220, 2004.

\bibitem{tartaruga2019experimental}
I.~Tartaruga, D.~Barton, D.~Rezgui, and S.~Neild.
\newblock Experimental bifurcation analysis of a wing profile.
\newblock In {\em International Forum on Aeroelasticity and Structural
  Dynamics: IFASD 2019}, 2019.

\end{thebibliography}

\end{document}